\def\IA{{\Bbb A}}

\def\IS{{\Bbb S}} 
\def\IK{{\Bbb K}}

\def\IC{\Bbb C} 
\def\ID{{\Bbb D}}

\def\zbar{{\overline{z}}} 
 
\def\wbar{{\overline{w}}}

\documentclass[12pt]{article} 

\usepackage[psamsfonts]{amssymb} 
\usepackage{amsfonts,amsmath}

\usepackage{epsfig,multicol}

\newtheorem{theorem}{Theorem}%[section] 
\newtheorem{lemma}{Lemma}%[section] 
\newtheorem{corollary}{Corollary}%[section] 

\newtheorem{conjecture}{Conjecture}[section]

\title{The Teichm\"uller problem for $L^p$-means of distortion.} 
\author{  Gaven J. Martin and Cong Yao} 
\date{} 
\begin{document}

\maketitle

\begin{abstract} Teichm\"uller's problem from 1944 is this: Given $x\in [0,1)$ find and describe the extremal quasiconformal map $f:\ID\to\ID$, $f|\partial \ID=identity$ and $f(0)=-x\leq 0$.  We consider this problem in the setting of minimisers of $L^p$-mean distortion. The classical result is that there is an extremal map of Teichm\"uller type with associated holomorphic quadratic differential having a pole of order one at $x$,  if $x\neq 0$.  For the $L^p$-norm, when $p=1$ it is known that there can be no locally quasiconformal minimiser unless $x=0$. Here we show that for $1\leq p<\infty$ there is a minimiser in a weak class and an associated Ahlfors-Hopf holomorphic quadratic differential with a pole of order $1$ at $f(0)=r$.  However,  this minimiser cannot be in $W^{1,2}_{loc}(\ID)$ unless $r=0$ and $f=identity$.  Hence there is no locally quasiconformal minimiser.  A similar statement holds for minimsers of the exponential norm of distortion.  We also use our earlier work to show that as $p\to\infty$,  the weak $L^p$-minimisers converge locally uniformly in $\ID$ to the extremal quasiconformal mapping,  and that as $p\to 1$  the weak $L^p$-minimisers converge locally uniformly in $\ID$ to the identity. \end{abstract}

\section{Introduction.}  Teichm\"uller's problem originated in his 1944 paper \cite{Teich}.  There is a nice commentary on this paper by V. Alberge \cite[Chapter 23]{Pap}.  The problem is this:

\medskip

\noindent{\bf Teichm\"uller's problem.}  Find and describe the quasiconformal map $f:\ID\to\ID$ from the unit disk $\ID$ to itself with
\begin{itemize}
\item  $f|\partial \ID=identity$
\item $f(0)=-x$,  $x\in(0,1)$.
\item $f$ has minimal distortion among all such mappings,
\end{itemize}

Teichm\"uller solved this problem by reducing it to the Gr\"otzsch problem for mappings between quadrilaterals after constructing certain double branched covers of $\ID$ over $0$ and over $-x$. It is important to note that the existence and the uniqueness of the extremal mapping cannot be deduced from the  Teichm\"uller theorem, \cite{Teich2}.  However,  the extremal mapping here is a Teichm\"uller mapping;
\begin{equation}
\mu_f(z) = \frac{f_\zbar}{f_z} = k \;\frac{\bar\phi}{\phi}, \quad \quad \mbox{ where $k\in [0,1)$ and $\phi$ is meromorphic.}
\end{equation}
In this case, $\phi$ has a pole of order $1$ at $0\in\ID$.
In fact in his 1962 ICM talk \cite{Jenkins} Jenkins wrote: ``Teichm\"uller enunciated the intuitive principle that the solution of a certain type of extremal problem for univalent  functions is determined by a quadratic differential for which the following prescriptions hold. If the competing mappings are to have a certain fixed point the quadratic differential will have a simple pole there.'' This will be an important observation in what follows. Generalising these ideas a little Kra defined a metric (now known as the Kra metric) on a Riemann surface $\Sigma$ as 
\begin{equation}
k_\Sigma(z_1,z_2) = \inf_{f}  \; \log \sqrt{K_f}, \quad\quad K_f = \Big\| \frac{1+|\mu_f(z)|}{1-|\mu_f(z)|} \Big\|_{L^\infty(\Sigma)}<\infty,
\end{equation}
where the infimum is taken over all $K_f$-quasiconformal self mappings $f:\Sigma\to\Sigma$ homotopic to the identify and such that $f(z_1)=z_2$.  That this is related to the Teichm\"uller problem can be seen by lifting $f$ to the universal cover $\hat{f}:\ID\to\ID$ where $\hat{f}|\partial \ID$ will be the identity if,  for instance, $\Sigma$ has finite area.  Further, the projections are locally conformal and therefore do not change the maximal distortion.   Modulus of continuity estimates for both $f$ and $f^{-1}$ (also quasiconformal) easily imply that the infimum is attained in the Kra metric.

\medskip We now state the problem we address in this article. A mapping $f:\ID\to \ID$ has {\em finite distortion} if 
\begin{enumerate}
\item $f\in W^{1,1}_{loc}(\ID)$,  the Sobolev space of functions with locally integrable first derivatives,
\item the Jacobian determinant $J(z,f)\in L^{1}_{loc}(\ID)$, and 
\item there is a measurable function ${\bf K}(z)\geq 1$, finite almost everywhere, such that 
 \begin{equation}\label{1.1}
 |Df(z)|^2 \leq {\bf K}(z) \, J(z,f), \hskip10pt \mbox{ almost everywhere in $\ID$}.
 \end{equation}
\end{enumerate}
See \cite[Chapter 20]{AIM} for the basic theory of mappings of finite distortion and the associated governing equations; degenerate elliptic Beltrami systems.
In (\ref{1.1}) the operator norm is used.  However this norm loses smoothness at crossings of the eigenvalues and for this reason when considering minimisers of distortion functionals one considers the distortion functional
\begin{equation}\label{1.2}
\IK(z,f) = \left\{\begin{array}{cc} \frac{\|Df(z)\|^2}{J(z,f)}, & \mbox{if $J(z,f)\neq 0$} \\
1, & \mbox{if $J(z,f)= 0$.} \end{array}\right.
\end{equation}
This was already realised by Ahlfors in his seminal work proving Teichm\"uller's theorem and establishing the basics of the theory of quasiconformal mappings, \cite[\S 3, pg 44]{Ahl}. We reconcile (\ref{1.1}) and (\ref{1.2}) by noting $\IK(z,f)=\frac{1}{2}\big({\bf K}(z)+1/{\bf K}(z) \big)$ almost everywhere,  where ${\bf K}(z)$ is chosen to be the smallest functions such that (\ref{1.1}) holds. 

\subsection{The $L^p$ Teichm\"uller Problem.}\label{Lp} Fix $p\geq 1$.  Find and describe the mapping  $f:\ID\to \ID$ of finite distortion attaining the infimum of
\begin{equation}\label{inf}
\inf _f \; \frac{1}{\pi} \int_\ID \IK(z,f)^p \; dz
\end{equation}
among all homeomorphisms of finite distortion with $f|\partial \ID=identity$ and $f(0)=-x$, $x\in[0,1)$.  

\medskip

We address this problem as follows.  First we establish existence: If $p> 1$,  any minimising sequence has a convergent subsequence in a regular a class of mappings described below.  We identify an associated meromorphic Ahlfors-Hopf quadratic differential and show it has a pole of order $1$.  Second we show that,  other than the identity,  none of these minimisers can be very regular - in particular they cannot be quasiconformal near $0$. Finally we show that as $p\to\infty$ these irregular solutions still converge in $W^{1,q}(\ID)$,  all $q<2$,  to the extremal quasiconformal mapping,  and that as $p\to 1$ these solutions must converge to the identity.   This is encapsulated in the following theorems.

\medskip

\noindent{\bf Data.} Let $f_0:\ID\to\ID$, $f|\partial\ID=identity$ and $f(0)=r\geq 0$ be the extremal quasiconformal mapping,  with distortion $K=\frac{1+k}{1-k}$.

\begin{theorem}\label{existence} Let $p>1$.  There is a mapping $f:\ID\to \ID$ of finite distortion attaining the infimum of
\begin{equation}\label{inf}
\inf _g \; \frac{1}{\pi} \int_\ID \IK(z,g)^p \; dz
\end{equation}
among all homeomorphisms $g:\ID\to\ID$ of finite distortion.  Further $f\in W^{1,q}(\ID,\ID)$ with the bound
\begin{equation}\label{5.2}
\|Df\|_{L^q(\ID)}\leq \pi  \Big(\frac{1+k^2}{1-k^2}\Big)^{\frac{p}{p+1}},
\end{equation} and there is $h:\overline{\ID}\to\overline{\ID}$ of finite distortion, which is
\begin{itemize} 
\item a continuous surjection from $\overline{\ID}\to \overline{\ID}$
\item $h|\partial \ID=identity$ and $h(r)=0$,  
\item $h\in W^{1,2}(\ID)$ is monotone (the preimage of a point is compact, connected and contractible)
\item The Ahlfors-Hopf quadratic differential 
\[ \Phi(w) = \IK(w,h)^{p-1} h_w \overline{h_\wbar} \]
is holomorphic in $\ID\setminus \{r\}$ and has a pole of order $1$ at $r$.
\item  There is a measurable set $X\subset\ID$ such that $|\ID-X|=0$, $(h\circ f)(z)=z$ for every $z\in X$, and $J(w,h)=0$ for almost every $w\in\ID-f(X)$.
\end{itemize}
\end{theorem}
We call $h$ a pseudo-inverse of $f$ and can therefore in a sense assert $f(0)=r$ and $f|\ID=identity$. We call $\Phi$ the Ahlfors-Hopf differential as these holomorphic quadratic differentials first appeared in Ahlfors seminal paper from 1953 \cite{Ahl} which gave the first complete proof of Teichm\"uller's theorem; he used the variation of $L^p$-means of distortion,  and then after a careful analysis let $p\to\infty$ !  

We expect uniqueness for minimisers in Theorem \ref{existence},  but have not yet established that for Ahlfors-Hopf differentials with a pole.
 
 \begin{theorem}\label{regularity}  Let $f$ be any minimiser given by Theorem \ref{existence}.  If either
 \begin{enumerate}
 \item $f\in W^{1,2}(U)$ for {\em any} open set $0\in U$, or
  \item $\IK(w,h) \in L^1(V)$ for {\em any} open set $z_0\in V$, 
 \end{enumerate}
 then $r=0$ and $f=h=identity$.
\end{theorem}

 \begin{theorem}\label{regularity}  Let $\{f_p\}_{1<p<\infty}$ be any family of minimisers given by Theorem \ref{existence}.   Then there are subsequences 
 \begin{itemize} 
 \item  $p_j\to\infty$ with $f_{p_j} \to f_\infty$ in $W^{1,r}(\ID)$, all $r<2$,  and $f_\infty$ is the extremal quasiconformal mapping of Teichm\"uller type, $f_\infty(0)=r$ and $h_{p_j} \to h_\infty$ in $W^{1,2}(\ID)$.
 \item   $q_j\to 1$ with $f_{q_j} \to identity$ in $W^{1,1}(\ID)$ and $h_{q_j} \to identity$ in $W^{1,2}(\ID)$.
 \end{itemize}
 \end{theorem}

In particular we have the following corollary.
\begin{corollary}
No minimiser for the $L^p$-Teichm\"uller problem can be locally quasiconformal other than the identity.
\end{corollary}

\noindent{\bf Remark.}. By way of contrast we make the following observations.
For the associated $L^p$-Gr\"otzsch problem,  $p>1$,  of minimising the $L^p$ mean distortion,
\begin{equation}
\frac{1}{\pi} \int_\IA \IK(z,f)^p \; dz, 
\end{equation} of homeomorphisms of finite distortion between annuli $f:\IA\to\IA'$, 
\[ \IA =\{z:1<|z|<R<\infty\},\quad \IA'=\{z:1<|z|<S<\infty\} \]
always has a minimiser.  Moreover that minimiser is a diffeomorphism of Teichm\"uller type,  meaning $\mu_f(z)=k(z) \frac{\bar\phi}{|\phi|}$,  $k(z)\in [0,1)$, $\phi$ holomorphic, \cite{MM}.  When $p=1$ a minimiser exists if and only if $R+1/R\leq 2S$.  This latter range of moduli is called the Nitsche range \cite{AIM}.  For the Gr\"otzsch problem of mappings between quadrilaterals,  that linear maps are the $L^p$ minimisers is,  more or less,  a direct consequence of the polyconvexity of the integrand $\IK(z,f)^p$,  $p\geq 1$. 

\medskip

In the case of the exponential means of distortion,  we have existence of minimisers because of the modulus of continuity of estimates \cite{AIMb}.  Yet these minimisers again cannot be locally quasiconformal.  Indeed they must be quite irregular.

\begin{theorem}\label{exp}  Fix $p\geq 0$, $x\in [0,1)$.  Then there is $f:\ID\to \ID$  with $f|\partial \ID=identity$ and $f(0)=-x$ a homeomorphism of finite distortion, minimising
\begin{equation}
\frac{1}{\pi} \int_\ID e^{p\IK(z,f)} \; dz
\end{equation}
among all homeomorphisms of finite distortion with $f|\partial \ID=identity$ and $f(0)=z_0$.  Further,  it cannot happen that $e^{q\IK(z,f)}\in L^1_{loc}(\ID)$ and $\IK(z,f) e^{p\IK(z,f)}\in L^1(\ID)$ for any $q>p$. \end{theorem}

We expect that the minimiser of Theorem \ref{exp} has a meromorphic Ahlfor's-Hopf differential with a pole of order $1$ at $0$, 
\begin{equation}\label{expphi}
\Phi(w) = e^{p\IK(w,h)} h_w(w) \overline{h_\wbar(w)},\quad\quad h=f^{-1}:\ID\to\ID.
\end{equation}
It is certainly the case that $\Phi$ is meromorphic in $\ID$ with a pole of order $1$ at $0$ if $e^{q\IK(z,f)}\in L^1_{loc}(\ID\setminus\{0\})$ for some $q>p$. A key issue is that it is not clear $\Phi$ as defined at (\ref{expphi}) lies in $L^1_{loc}(\ID)$.

\section{Existence.}
 We first discuss the natural class where one might find a minimiser. If $\{f_j\}_{j=1}^{\infty}$ is a minimising sequence of homeomorphic mappings of finite distortion, then H\"older's inequality gives the {\em a priori} bounds (see \cite{IMO})
\begin{eqnarray} \nonumber
\Big[\int_\ID \|Df_j(z)\|^{\frac{2p}{p+1}}dz\Big]^{p+1} & \leq & \int_{\ID }\IK^p(z,f_j)\; dz\; \cdot\; \Big[\int_\ID J(z,f_j)dz\Big]^p \\ & \leq & \pi^p\int_{\ID }\IK^p(z,f_j)\; dz. \label{1.7}
\end{eqnarray}
Thus there is a subsequence $f_j\rightharpoonup f$ weakly in $W^{1,\frac{2p}{p+1}}(\ID)$. Meanwhile, the sequence of inverses,  $h_j=f_j^{-1}$, satisfies
\begin{align} 
\int_\ID\|Dh_j(w)\|^2dw&=\int_\ID\IK(w,h_j)J(w,h_j)dw\notag\\
&\leq\int_{\ID }\IK^p(w,h_j)J(w,h_j)dw=\int_{\ID }\IK^p(z,f_j)dz,\label{1.8}
\end{align}
so they converge in $W^{1,2}(\ID)$. Note the change of variables formula follows from \cite{H,KM}. Such a sequence $h_j$ converges to a continuous function $h$ locally uniformly in $\ID$ \cite{GV,IM} dues to a modulus of continuity estimate for $W^{1,2}(\ID)$ mappings of finite distortion,  \cite{AIMb}.  In fact in \cite{IKO} it is proved that $h$ will be locally Lipschitz. However, on the $f$ side, functions in $W^{1,\frac{2p}{p+1}}(\ID)$ are not usually continuous. The remainder of Theorem \ref{existence},  obtaining the set $X$ of full measure and the vanishing of the Jacobian on this set are basically local and exactly follow the arguments of our earlier work \cite[\S 5]{MY}.

\section{Equations.}
Set $\ID^*=\ID\setminus\{0\} $.  It is essentially proved in \cite{IMO} that a minimiser of Problem \ref{Lp}  must satisfy the following inner-variational equation:
\begin{equation}\label{1.5}
2p\int_\ID \IK(z,f)^p\frac{\overline{\mu_f}}{1+|\mu_f|^2}\varphi_\zbar dz=\int_\ID \IK(z,f)^p\varphi_zdz,\quad\forall\varphi\in C_0^\infty(\ID^*).
\end{equation}
The only difference in that case was $\ID$ replacing $\ID^*$ but the calculation is local.  This equation arises as follows.  Let $\varphi\in C^{\infty}_{0}(\ID^*)$ with $\|\nabla\varphi\|_{L^{\infty}(\ID)}<1$.  Then for $t\in (-\frac{1}{2},\frac{1}{2})$ the mapping $g^t(z)=z+t\varphi(z)$ is a diffeomorphism of $\ID$ to itself which extends to the identify on the boundary $\IS$ and has $g^t(0)=0$.  If $f$ is a mapping of finite $L^p$-mean distortion distortion,  then so is $f\circ (g^t)^{-1}$ and $f\circ (g^t)^{-1}\in W^{1,2p/(p+1)}_0(\ID)$. If $f$ is a homeomorphism,  so is $f\circ (g^t)^{-1}$,  and if $f(0)=z_0$ then $f\circ (g^t)^{-1}(0)=z_0$.  Thus,  for each $t$,  $f\circ (g^t)^{-1}$ is a candidate for the extremal problem.

The function 
\begin{equation} t\mapsto \frac{1}{\pi} \int_\ID \IK(z, (f\circ (g^t)^{-1})^p\; dz\end{equation}  is a smooth function of $t$. Thus if  $f$ is a minimiser in any reasonable class (that is we may relax the assumption that $f$ is a homeomorphism) we have
\[ \frac{d}{dt} \Big|_{t=0} \frac{1}{\pi} \int_\ID \IK(z, (f\circ g^t))^p\; dz  = 0. \]
It is a nice calculation to verify that this equation is equivalent to (\ref{1.5}). As a hint, use the change of variables $z=g^t(w)$ and then the composition formula for Beltrami coefficients \cite[Theorem 5.5.6]{AIMb}. It is interesting to note that (\ref{1.5}) implies that $\mu_f$ is constant on any open set that $|\mu_f|$ is constant.  

Next, (\ref{1.5}) tells us that if we solve $F_\zbar=\IK(z,f)^p-1$,  which we can do using the Cauchy transform on $\ID$,  then we have
\begin{equation}\label{1.5}
 \int_\ID F_z \varphi_\zbar dz=\int_\ID F_\zbar \varphi_z dz,\quad\forall\varphi\in C_0^\infty(\ID^*).
\end{equation}
where
\[ F_z=\IK(z,f)^p\frac{2p \overline{ \mu_f}}{1+|\mu_f|^2} \]
and $F\in W^{1,1}(\ID)$.
This calculation is carried out carefully in \cite[Lemma 2.1]{MY}.  As there,  it leads to the following autonomous equation for $F$,
\begin{equation}\label{Feqn}
F_\zbar = A_p(|F_z|), \quad\quad 0\leq A_p'(t) \leq \frac{1}{p}.
\end{equation}
In \cite{HM} it is shown that $W^{1,s}_{loc}(\ID^*)$ solutions to this equation, $s>1+1/p$ are $C^{\infty}$ and both $F$ and $F_z$ are $ \frac{p+1}{p-1}$--quasiregular.  In fact since we do have $F$,  as observed in \cite{MY2} the injectivity of Beltrami operators in the borderline case \cite{AIS} and see \cite[\S 14.4.2]{AIMb} implies $F$ is smooth if it is in $W^{1,1+1/p}_{loc}(\ID^*)$.

\begin{lemma} \label{lemma1} If $U$ is open,  $0\in U$, and $\IK(z,f)^p \in L^{1+1/p}(U\setminus \{0\})$,  then $f|U$ is a diffeomorphism.
\end{lemma}
\noindent{\bf Proof.} We have already observed that in these circumstances $F$ is smooth and quasiregular as it implies $F_z\in L^{1+1/p}(U)$.  The Stoilow factorisation theorem \cite[\S5.5]{AIMb} tells us we can write $F=\phi\circ g$ where $g:\ID^*\to\ID^*$ is $ \frac{p+1}{p-1}$--quasconformal and $\phi:\ID^*\to\IC$ is holomorphic. The point $0$ is an isolated singularity for the quasiconformal map $g$ and is hence removable and $g\in W^{1,1+p}(\ID)$ -- again using borderline injectivity. Then for $z\in \ID^*$
\begin{eqnarray*} \phi'(z) &=&  (F\circ g^{-1})_z = F_z(g^{-1})(g^{-1})_w+F_\zbar(g^{-1})\overline{(g^{-1})_\wbar} \\
 \int_\ID |F_z(g^{-1})(g^{-1})_w dw| & = &  \int_\ID |F_z g_z dz|   \leq \Big(\int_\ID |F_z|^{1+1/p} dz\Big)^{\frac{p}{p+1}} \Big(\int_\ID |g_z|^{1+p} dz   \Big)^{\frac{1}{p+1}} \\
 & \leq  &  C \Big(\int_\ID \IK(z,f)^{p+1} dz\Big)^\frac{p}{p+1}, \quad\quad C>>1.
 \end{eqnarray*}
Hence $\phi'\in L^{1}(U)$ and holomorphic on $U\setminus\{0\}$, and therefore $\phi$ has $0$ as a removable singularity.  Thus  $F$ is smooth and bounded in $U$ as $F$ also solves (\ref{Feqn}). Hence $\mu_f$ is smooth and $f|U$ is a diffeomorphism.  \hfill $\Box$

\medskip

We now argue similarly on the other side.

\begin{lemma}\label{lemma2} If $V$ is open,  $r \in V$, and $\IK(w,h) \in L^{1}(V)$,  then $f|U$ is a diffeomorphism. 
\end{lemma}
\noindent{\bf Proof.} It was essentially first proved by Ahlfors \cite{Ahl}  that the function $h$ a pseudo inverse of minimiser of Problem \ref{Lp}  must satisfy the following inner-variational equation:
\begin{equation}\label{1.6}
\int_\ID \IK^{p-1} h_w \overline{h_\wbar} \varphi_\zbar dz, \quad\forall\varphi\in C_0^\infty(\ID^*).
\end{equation}
Since
\[ \IK^{p-1}(w,h) h_w \overline{h_\wbar} = \IK^{p-1}(w,h) |h_w|^2  \overline{\mu_h} =  \IK^{p}(w,h)J(w,h) \frac{\overline{\mu_h}}{1+|\mu|^2} \in L^1(\ID), \]
the Weyl lemma implies the Ahlfors-Hopf differential
\begin{equation}\label{AH} \Phi(w) = \IK^{p-1}(w,h) h_w \overline{h_\wbar} \end{equation}
is holomorphic in $\ID^*$ with a pole of at most order $1$ at $r$.

This equation arises just as above.  Let $\varphi\in C^{\infty}_{0}(\ID\setminus \{r\})$ with $\|\nabla\varphi\|_{L^{\infty}(\ID)}<1$.  Then for $t\in (-\frac{1}{2},\frac{1}{2})$ the mapping $g^t(z)=z+t\varphi(z)$ is a diffeomorphism of $\ID$ to itself which extends to the identify on the boundary $\IS$ and has $g^t(r)=r$.  For each $t$,  $h\circ (g^t)^{-1}$ is a candidate for the extremal problem. The function 
\begin{equation} t\mapsto \frac{1}{\pi} \int_\ID \IK(z, (h\circ (g^t)^{-1})^p\;J(z,h\circ (g^t)^{-1})  dz \end{equation}  is a smooth function of $t$ and 
\[ \frac{d}{dt} \Big|_{t=0} \frac{1}{\pi} \int_\ID \IK(z, (h\circ g^t))^p\; J(z,f\circ (g^t)^{-1}) \; dz  = 0, \]
leads to (\ref{1.6}). Now $h\in W^{1,2}(\ID)$ is a monotone mapping of finite distortion with $\IK(w,h)\in L^1(V)$,  by hypothesis. Thus \cite{OZ} implies $h$ is discrete and open,  and hence $h|V$ is a homeomorphism.  Set $h(V)=U$ and note that if $|\mu_h|\geq \frac{3-\sqrt{5}}{2}$,  then
\[ \IK^{p+1}(w,h)J(w,h) = \IK(w,h) \Phi(w) \frac{1+|\mu|^2}{|\mu_h|} \leq 3 \IK(w,h) \Phi(w) \in L^1(V)\]
since $\Phi$ is continuous in $\ID$.  On the other hand,  if $|\mu_h|\leq \frac{3-\sqrt{5}}{2}$,  then 
\[ \IK(w,h)^{p+1}J(w,h) \leq \frac{ 3}{\sqrt{5}}J(w,h) \in L^1(\ID) \]
 We thus have $ \IK(w,h)^{p+1}J(w,h)\in L^1(V)$ and hence
\[ \int_V\IK(w,h)^{p+1}J(w,h) \;dw = \int_U\IK(z,f)^{p+1}\; dz, \] 
and so Lemma \ref{lemma1} now applies to show $f$ is a diffeomorphism.  \hfill $\Box$

\section{Proofs and results.}

\subsection{Proof of Theorem \ref{regularity} .} We have shown above in Lemmas \ref{lemma1},\ref{lemma2} that under either of the hypotheses of Theorem \ref{regularity} the mapping $f$ is a diffeomorphism near $0$. Hence $h$ is a diffeomorphism near $r$ and hence the Ahlfors-Hopf differential $\Phi$ is holomorphic in $\ID$.  However,  there is a unique monotone mapping in $W^{1,2}(\ID)$ with holomorphic Ahlfors-Hopf differential in $\ID$ and $h|\partial \ID=identity$,  see \cite{MY2}.  It is of course the identity.  Thus $r=0$ and $f=identity$.

\subsection{Limiting regimes.}  The arguments concerning what happens as $p\to\infty$ in large part go back to Ahlfors again.  The case $p=1$ is also by now well-known \cite{AIMO}.  But the direct methods suggested here apparently fail for we do not get a uniform elliptic estimate and we must resort to an alternative approach. In \cite{MY2} we prove that as $p\to 1$ the psuedo-inverses $h_p$ of minimisers $f_p$ converge locally uniformly in $\ID\setminus \{r\}$ to a harmonic mapping $h_1$ in $\ID\setminus \{r\}$.  The limit $h_1$ is continuous due to the modulus of continuity estimate independent of $p\geq 1$.  Then $r$ is a removable singularity and $h_1=identity$.  It directly follows that the convergence is locally uniform in $\ID$ as stated (consider a small circle around $r$ and the modulus of continuity). When $p\to\infty$ we show the local uniform limit exists and is an extremal quasiconformal mapping and further we identify when the approximating sequence is a Hamilton sequence,  making this limit a uniquely extremal Teichm\"uller mapping for its boundary values.These are carefully spelled out in \cite[\S 7]{MY} and we leave the reader to explore that section.

 Second,  the only uniqueness statement we currently know requires $\Phi\in L^1(\ID)$.  However,  if $e^{q\IK(w,h)} \in L^1_{loc}(\ID^*)$ for any $q>p$ and $\Phi\in L^1(\ID)$,  then both these problems can be resolved and $f$ (and hence $h$) will be a diffeomorphism, \cite{MY4} and the proofs follow the same argument.

\section{The Kra metric.} From what we have above,  it follows that given $z_0,z_1\in \ID$  the problem describing the mapping  $f:\ID\to \ID$ of finite distortion attaining the infimum of
$
\inf _f \; \frac{1}{\pi} \int_\ID \IK(z,f)^p \; dz$
among all homeomorphisms of finite distortion homotopic to the identity and $f(z_0)=z_1$ has exactly the same outcomes.  

 \subsection{$L^p$ extremals on a Riemann surface.} 
 For a Riemann surface  $\Sigma$ of finite area and $z_0,z_1\in \Sigma$ the problem is a little more interesting when minimising
\[   \inf _f \;   \int_\Sigma \IK(z,f)^p \; d\sigma(z).  \]   
As above, modulus of continuity arguments show that there is a monotone $h:\Sigma\to\Sigma$ in $W^{1,2}(\Sigma,\Sigma)$ with meromorphic Ahlfors-Hopf differential with a pole of order $1$ at $z_0$.  If $\IK(w,h)\in L^1(\ID)$,  then $h$ is a diffeomorphism and there is no pole.  To see this simply argue locally as above.  Thus $f=h^{-1}:\Sigma\to\Sigma$ is a critical point for the problem of identifying the extremal mapping for the $L^p(\Sigma)$-mean distortion (in the homotopy class of the identity).  We expect uniqueness for this problem.  Note that now $f$ and $h$ lift to diffeomorphisms of $\ID$ which are the identity on the boundary.   The lift, say $\hat{h}$, also has a holomorphic Ahlfors-Hopf differential $\hat{\Phi}$.

\begin{conjecture}\label{conj} Suppose $\hat{h}:\ID\to\ID$ is a diffeomorphism of $\ID$, a homeomorphism of $\overline{\ID}$ and $\hat{h}|\partial\ID=identity$.  If $\hat{\Phi} = \IK(w,h)^{p-1} h_w \overline{h_\wbar}\, \eta(h)$ is holomorphic,  then $\hat{h}=identity$.
\end{conjecture}
Here $\eta$ is the hyperbolic area density on $\Sigma$.

We know the conjecture is true in the case $\hat{\Phi}\in L^1(\ID)$,  but in this case we see that if $P$ is a convex hyperbolic fundamental domain for the action of the fundamental group $\Gamma$ of $\Sigma$ by hyperbolic isometries of $\ID$ we may calculate that
\begin{eqnarray*} \int_\ID |\hat{\Phi}| & = &  \int_{\cup_{\gamma\in \Gamma}\gamma(P)} |\hat{\Phi}| =\sum_{\gamma\in \Gamma}\int_{ \gamma(P)} |\hat{\Phi}| \\ & =& \sum_{\gamma\in \Gamma}\int_{P} |\hat{\Phi}(\gamma)||\gamma'|^2 =   \sum_{\gamma\in \Gamma}\int_{P} |\hat{\Phi}| = +\infty,
\end{eqnarray*} 
unless $\Phi=0$. In a later article we expect to address this problem of uniqueness on Riemann surfaces of finite area.

 \subsection{The exponential extremals on a Riemann surface. }  We have noted there is a minimiser because of modulus of continuity estimates for both $f$ and its inverse $h$, \cite{AIMb}. These equicontinuity estimates carry over to mappings of finite distortion between Riemann surfaces. One may formally perform the variational calculations as above to find the pointwise defined differential
 \[ \Phi(w) = e^{p\IK(z,h)} h_w \overline{h_\wbar}\eta(h).   \] 
There are two issues.  First $f$ might be non-variational, \cite{MY3}.  This means that  $e^{p\IK(z,f)}\IK(z,f) \not\in L^{1}_{loc}(\Sigma)$,  so 
\[ \Phi(w) = e^{p\IK(z,h)} h_w \overline{h_\wbar}\eta(h)  \sim \; \IK(z,h) e^{p\IK(w,h)} J(w,h)  \not\in L^1_{loc}(\Sigma) \]
and so we cannot directly conclude it is holomorphic.  However,  in the setting of Riemann surfaces this is addressed in  \cite[Theorem 3]{MY5} using the Riemann-Roch theorem and only minor changes are necessary to those arguments to conclude that the the case at hand $\Phi$ is indeed a meromorphic  Ahlfors-Hopf differential with a pole of order $1$ at $z_0$ and holomorphic in $\Sigma\setminus\{z_0\}$.  Then,  as before regularity at $z_0$ implies that $h$ is a diffeomorphism. Precisely we would require $e^{q\IK(z,f)}\in L^1(U)$ for some open $U$ containing $z_0$ and any $q>p$.  Thus we return to the problem of uniqueness as per Conjecture \ref{conj}.

\end{document}